\documentclass[12pt,letterpaper]{article}

\usepackage{latexsym} 
\usepackage{array}
\usepackage{enumerate}
\usepackage{theorem}

\newtheorem{theorem}{Theorem}

\begin{document}

\title{Creator-annihilator domains \\ and the number operator}
\author{P. L. Robinson}
\date{}
\maketitle

\section*{Abstract} 

We show that for the bosonic Fock representation in infinite dimensions, the maximal common domain of all creators and annihilators properly contains the domain of the square-root of the number operator. 

\section*{Introduction}

A standard construction of the bosonic Fock representation of a complex Hilbert space $V$ is in terms of creators and annihilators on the symmetric algebra $SV$ in which the number operator $N$ scales elements by homogeneous degree; the symmetric algebra is completed relative to a canonical inner product and these various operators are extended to maximal domains in the resulting Fock space $S[V]$. 

\medbreak

It is well-known that all creators and annihilators are defined on the domain of the square-root $N^{\frac{1}{2}}$. A natural question (raised on page 16 of [2] by Berezin, among others; compare page 65 of [4]) is whether the domain of $N^{\frac{1}{2}}$ coincides with the maximal common domain of all creators and annihilators. Here, we demonstrate that the answer to this question is affirmative when $V$ is finite-dimensional but negative when $V$ is infinite-dimensional. 

\medbreak 

As noted above, the specific question that is answered in this paper appears in [2]. Among many standard references concerning the bosonic Fock representation, we cite [1] and [3]. The particular approach taken here (involving the full antidual of the symmetric algebra) was introduced in [6] as a means to establishing a generalized version of the classic Shale theorem on the implementation of symplectic automorphisms, independently of the generalization previously presented in [5]. 

\section*{Fock-lore}

We begin by recalling certain familiar elements of Fock-lore, pertaining to the construction of bosonic Fock space and the various operators defined therein. For traditional accounts, see a standard text such as [1] or [3]; for an account in line with the present paper, see [6]. 

\medbreak 

To be explicit, let $V$ be a complex Hilbert space. Extend its inner product $\langle \cdot \vert \cdot \rangle$ to the symmetric algebra $SV = \bigoplus_{n \geq 0} S^n V$ by declaring that the homogeneous summands $(S^n V : n \geq 0)$ be perpendicular and that if $n \geq 0$ and $x_1, \dots , x_n, y_1, \dots , y_n \in V$ then 
\nonumber 
\[
\langle x_1 \cdots x_n \vert y_1 \cdots y_n \rangle = \sum_{p} \prod_{j = 1}^{n} \langle x_j \vert y_{p(j)} \rangle 
\]
where $p$ runs over all permutations of $\{1, \dots , n \}$; the resulting complex Hilbert space completion is (by definition) the bosonic Fock space $S[V] = \bigoplus_{n \geq 0} S^n [V]$. 

\medbreak 

For our purposes, it is convenient to introduce also the full antidual $SV'$ comprising all antilinear functionals $SV \rightarrow \mathbf{C} $. This full antidual $SV'$ is naturally a commutative, associative algebra under the product defined by 
\nonumber 
\[
\Phi, \Psi \in SV' \Longrightarrow [\Phi \Psi](\theta) = (\Phi \otimes \Psi)(\Delta \theta) 
\]
where the coproduct $\Delta : SV \rightarrow SV \otimes SV$ arises when the canonical isomorphism $S(V \oplus V) \equiv SV \otimes SV$ follows the homomorphism $SV \rightarrow S(V \oplus V)$ induced by the diagonal map $V \rightarrow V \oplus V$. The inner product on $SV$ engenders an algebra embedding 
\nonumber 
\[
SV \rightarrow SV': \phi \mapsto \langle \cdot \vert \phi \rangle
\]
and $S[V]$ is identified with the subspace of $SV'$ comprising all bounded anti-linear functionals. 

\medbreak 

Let $v \in V$. The creator $c(v)$ is defined initially on $SV$ as the operator of multiplication by $v$: 
\nonumber 
\[
\phi \in SV \Longrightarrow c(v) \phi = v \phi.
\]
The annihilator $a(v)$ is defined initially on $SV$ as the linear derivation that kills the vacuum $\textbf{1} \in \textbf{C} = S^0 V$ and sends $w \in V = S^1 V$ to $\langle v \vert w \rangle$: thus, if $v_1, \dots , v_n \in V$ then 
\nonumber 
\[
a(v) (v_1 \cdots v_n) = \sum_{j = 1}^{n} \langle v \vert v_j \rangle v_1 \cdots \widehat{v_j} \cdots v_n 
\]
where the circumflex $\: \widehat{\cdot}\: $ signifies omission. As may be verified by direct calculation, $c(v)$ and $a(v)$ are mutually adjoint: 
\nonumber 
\[
\phi, \psi \in SV \Longrightarrow \langle \psi \vert c(v) \phi \rangle = \langle a(v) \psi \vert \phi \rangle .
\]
Accordingly, the creator $c(v)$ and annihilator $a(v)$ extend to $SV'$ by antiduality: if $\Phi \in SV'$ and $\psi \in SV$ then 
\nonumber 
\[
[c(v) \Phi](\psi) = \Phi [a(v) \psi]
\]
and 
\nonumber 
\[
[a(v) \Phi](\psi) = \Phi [c(v) \psi].
\]
Finally, the corresponding (mutually adjoint) operators in Fock space $S[V] \subset SV'$ are defined by restriction to the (coincident) natural domains 
\nonumber 
\[
\mathcal{D}[c(v)] = \{ \Phi \in S[V] : c(v) \Phi \in S[V] \}
\]
and 
\nonumber 
\[
\mathcal{D}[a(v)] = \{ \Phi \in S[V] : a(v) \Phi \in S[V] \}.
\]

\medbreak 

The number operator $N$ is defined initially on $SV$ by the rule 
\nonumber 
\[
n \geq 0 \Longrightarrow N \vert S^n V = n I
\]
and extends to $SV'$ by antiduality: 
\nonumber 
\[
\Phi \in SV' , \psi \in SV \Longrightarrow [N \Phi] (\psi) = \Phi [N \psi].
\]
The number operator in $S[V] \subset SV'$ is defined by restriction to the natural domain 
\nonumber 
\[
\mathcal{D}[N] = \{\Phi \in S[V] : N \Phi \in S[V] \} 
\]
which may be identified in terms of the decomposition $S[V] = \bigoplus_{n \geq 0} S^n [V]$ as 
\nonumber 
\[
\mathcal{D}[N] = \{\sum_{n \geq 0} \Phi_n \in S[V] : \sum_{n \geq 0} \Vert n \Phi_n \Vert^2 < \infty \}.
\]
We remark that the number operator $N$ in $S[V]$ is selfadjoint (indeed, positive): $\bigoplus_{n \geq 0} S^n [V]$ is its spectral decomposition, whence powers of $N$ are readily described in concrete terms; in particular, 
\nonumber 
\[
\mathcal{D}[N^{\frac{1}{2}}] = \{ \Phi \in S[V] : \sum_{n \geq 0} n \Vert \Phi_n \Vert^2 < \infty \}. 
\]

\section*{Theorems}

Having established sufficient background, we now proceed to our primary task: that of relating $\mathcal{D}[N^{\frac{1}{2}}]$ to the maximal common domain of all creators and annihilators in Fock space. 

\medbreak 

Let $u \in V$ be (for convenience) a unit vector: the unitary decomposition $V = \mathbf{C}u \bigoplus u^{\perp}$ induces (for each $n \geq 0$) a unitary decomposition 
\nonumber 
\[
S^n V = \bigoplus_{p + q = n} \{S^p (\mathbf{C}u) \otimes S^q (u^{\perp}) \}.
\]
Decomposing $\phi \in S^n V$ as 
\nonumber 
\[
\phi = \sum_{p + q = n} u^p \otimes \psi_q 
\]
we note (by perpendicularity) that 
\nonumber 
\[
\Vert \phi \Vert^2 = \sum_{p + q = n} \Vert u^p \Vert^2 \Vert \psi_q \Vert^2 = \sum_{p + q = n} p! \: \Vert \psi_q \Vert^2 
\]
and that 
\nonumber 
\[
\Vert c(u
) \phi \Vert^2 = \sum_{p + q = n} \Vert u^{p + 1} \Vert^2 \Vert \psi_q \Vert^2 = \sum_{p + q = n} (p + 1)! \: \Vert \psi_q \Vert^2 
\]
whence it follows that 
\nonumber 
\[
\Vert c(u) \phi \Vert^2 \leq (n + 1) \Vert \phi \Vert^2.
\]
This inequality continues to apply when $\phi$ lies in the closure $S^n [V]$; in particular, $S^n [V] \subset \mathcal{D}[c(u)]$. Now, if $\Phi \in \mathcal{D}[N^\frac{1}{2}]$ then 
\nonumber 
\[
\Vert c(u) \Phi \Vert^2 = \sum_{n \geq 0}\Vert c(u) \Phi_n \Vert^2 \leq  \sum_{n \geq 0} (n + 1) \Vert \Phi_n \Vert^2 = \Vert N^\frac{1}{2} \Phi \Vert^2 + \Vert \Phi \Vert^2 
\]
whence $\Phi \in \mathcal{D}[c(u)] = \mathcal{D}[a(u)]$. Lifting the convenient hypothesis that $u \in V$ be a unit vector, we have justified the following result.

\begin{theorem}
If $u \in V$ then $\mathcal{D}[N^\frac{1}{2}]$ is contained in $\mathcal{D}[c(u)] = \mathcal{D}[a(u)]$. 
\end{theorem}

\medbreak 

Thus, $\mathcal{D}[N^\frac{1}{2}]$ is contained in the maximal common domain of all creators and annihilators in Fock space. 

\medbreak 

We now consider the reverse containment in case $V$ is finite-dimensional, with $(u_1, \dots , u_m)$ a unitary basis. If the integers $n_1, \dots , n_m \geq 0$ have sum $n$ then for each $j \in \{ n_1, \dots , n_m \}$ 
\nonumber 
\[
c(u_j) a(u_j) (u_{1}^{n_1} \cdots u_{m}^{n_m}) = n_j (u_{1}^{n_1} \cdots u_{m}^{n_m})
\]
so  
\nonumber 
\[
\sum_{j = 1}^{m} c(u_j) a(u_j) (u_{1}^{n_1} \cdots u_{m}^{n_m}) = n (u_{1}^{n_1} \cdots u_{m}^{n_m}).
\]
By linearity, it follows that 
\nonumber 
\[
\sum_{j = 1}^{m} c(u_j) a(u_j)\:  \phi = n \phi
\]
whenever $\phi \in S^n V$ and indeed whenever $\phi \in S^n [V]$ by the boundedness of creators and annihilators on homogeneous elements of Fock space. Now, if $\Phi \in \mathcal{D}[a(u_1)] \cap \cdots \cap \mathcal{D}[a(u_m)]$ then 
\nonumber 
\[
\infty > \sum_{j = 1}^{m} \Vert a(u_j) \Phi \Vert^2 = \sum_{j = 1}^{m} \sum_{n \geq 0} \Vert a(u_j) \Phi_n \Vert^2 
\]
whence (valid) passage of annihilators across the inner product as creators yields 
\nonumber 
\[
\infty > \sum_{n \geq 0}\langle \Phi_n \vert \sum_{j = 1}^{m} c(u_j) a(u_j) \Phi_n \rangle = \sum_{n \geq 0} n \Vert \Phi_n \Vert^2
\]
which places $\Phi$ in $\mathcal{D}[N^\frac{1}{2}]$. This justifies the following result.

\begin{theorem} 
If $(u_1, \dots , u_m)$ is a unitary basis for the finite-dimensional $V$ then 
\nonumber 
\[
\mathcal{D}[N^\frac{1}{2}] = \mathcal{D}[a(u_1)] \cap \cdots \cap \mathcal{D}[a(u_m)]
\]
and if $\Phi$ lies in this domain then  
\nonumber 
\[
\Vert N^\frac{1}{2} \Phi \Vert^2 = \sum_{j = 1}^{m} \Vert a(u_j) \Phi \Vert^2.
\]

\end{theorem}

\medbreak 

In particular, if $V$ is finite-dimensional then $\mathcal{D}[N^\frac{1}{2}]$ coincides with the maximal common domain of all creators and annihilators in Fock space. 

\medbreak 

The case in which $V$ is infinite-dimensional is different. To see this, let $\Phi = \sum_{n > 0} \Phi_n \in SV'$ be defined by 
\nonumber 
\[
n > 0 \Longrightarrow \Phi_n = \lambda_n (u_n)^n / n! 
\]
where $(\lambda_n : n > 0)$ is a complex sequence and where the unit vectors $( u_n : n > 0 )$ in $V$ are perpendicular. From 
\nonumber 
\[
n > 0 \Longrightarrow \Vert \Phi_n \Vert^2 = \vert \lambda_n \vert^2 \Vert (u_n)^n \Vert^2 / (n!)^2 = \vert \lambda_n \vert^2 / n! 
\]
it follows that 
\nonumber 
\[
\Phi \in S[V] \Longleftrightarrow \sum_{n > 0} \frac{\vert \lambda_n \vert^2}{n!} < \infty 
\]
and that 
\nonumber 
\[
\Phi \in \mathcal{D}[N^\frac{1}{2}] \Longleftrightarrow \sum_{n > 0} \frac{\vert \lambda_n \vert^2}{(n - 1)!}< \infty.
\]
Further, let $v \in V$: if $n > 0$ then 
\nonumber 
\[
a(v) \Phi_n = \langle v \vert u_n \rangle \lambda_n (u_n)^{n - 1}/(n - 1)!
\]
whence 
\nonumber 
\[
\Vert a(v) \Phi \Vert^2 = \sum_{n > 0} \vert \langle v \vert u_n \rangle \vert^2 \frac{\vert \lambda_n \vert^2}{(n - 1)!}.
\]
Now, if $(\lambda_n : n > 0)$ is chosen so that 
\nonumber 
\[
\sum_{n > 0} \vert \lambda_n \vert^2 / n! < \infty = \sum_{n > 0} \vert \lambda_n \vert^2 / (n - 1)!
\]
then $\Phi \in S[V] \setminus \mathcal{D}[N^\frac{1}{2}]$ while if also $(\vert \lambda_n \vert^2 / (n - 1)! : n > 0)$ is bounded above by $K > 0$ then 
\nonumber 
\[
\Vert a(v) \Phi \Vert^2 \leq \sum_{n > 0} \vert \langle v \vert u_n \rangle \vert^2 K \leq K \Vert v \Vert^2  
\]
which places $\Phi$ in $\mathcal{D}[a(v)] = \mathcal{D}[c(v)]$. Of course, these conditions are easily satisfied: for example, when  $n > 0$ simply take $\vert \lambda_n \vert^2 / (n - 1)! = 1/n$. As the vector $v \in V$ is arbitrary, the following result is justified. 

\begin{theorem} 
If $V$ is infinite-dimensional then $\mathcal{D}[N^\frac{1}{2}]$ is properly contained in the maximal common domain of all creators and annihilators. 
\end{theorem} 

\medbreak 

Without proof, we remark that similar results hold for higher powers of the number operator: thus, if $k > 0$ then $\mathcal{D}[N^\frac{k}{2}]$ is contained in the maximal common domain of all degree $k$ polynomials in creators and annihilators, containment being strict precisely when $V$ is infinite-dimensional.

\section*{References}

\noindent
[1] J. C. Baez, I. E. Segal and Z. Zhou, \textit{Introduction to Algebraic and Constructive Quantum Field Theory}. Princeton University Press (1992). 

\medbreak 

\noindent
[2] F. A. Berezin, \textit{The Method of Second Quantization}. Academic Press (1966). 

\medbreak 

\noindent 
[3] O. Bratteli and D. W. Robinson, \textit{Operator Algebras and Quantum Statistical Mechanics II}. Springer-Verlag (1981). 

\medbreak 

\noindent 
[4] J. T. Ottesen, \textit{Infinite Dimensional Groups and Algebras in Quantum Physics}. Springer-Verlag (1995). 

\medbreak 

\noindent 
[5] S. M. Paneitz, J. Pedersen, I. E. Segal and Z. Zhou, \textit{Singular Operators on Boson Fields as Forms on Spaces of Entire Functions on Hilbert Space}. J. Functional Analysis \textbf{100} (1991) 36-58. 

\medbreak 

\noindent
[6] P. L. Robinson, \textit{The bosonic Fock representation and a generalized Shale theorem}. University of Florida preprint (1998); arXiv:1203.5841v1. 

\bigbreak 

\begin{flushleft}

Department of Mathematics 

University of Florida 

Gainesville, FL 32611 

\bigbreak 

e-mail: paulr@ufl.edu

\end{flushleft}

\end{document}